\documentclass[11pt]{article}
\usepackage{amscd}
\newcommand{\al}{\alpha}
\newcommand{\be}{\beta}
\newcommand{\ga}{\gamma}
\newcommand{\de}{\delta}
\newcommand{\ep}{\epsilon}
\newcommand{\dbar}{\overline{\partial}}
\newcommand{\dd}[1]{\partial_{#1}}
\newcommand{\dbr}[1]{\partial_{\overline{#1}}}

\newcommand{\ddt}[1]{\frac{\partial #1}{\partial t}}
\newcommand{\R}[4]{R_{#1 \overline{#2} #3 \overline{#4}}}
\newcommand{\Ric}[2]{R_{#1 \overline{#2}}}

\newcommand{\lt}{\tilde{\triangle}}
\newcommand{\chii}[2]{\chi^{#1 \overline{#2}}}
\newcommand{\ch}[2]{\chi_{#1 \overline{#2}}}
\newcommand{\chz}[2]{\chi_{0 \, #1 \overline{#2}}}
\newcommand{\g}[2]{g_{#1 \overline{#2}}}
\newcommand{\gi}[2]{g^{#1 \overline{#2}}}
\newcommand{\h}[2]{h^{#1 \overline{#2}}}

\begin{document}
\title{Convergence of the J-flow on Kahler surfaces}
\newcounter{theor}
\setcounter{theor}{1}
\newtheorem{theorem}{Theorem}[section]
\newtheorem{proposition}{Proposition}[section]
\newtheorem{lemma}{Lemma}[section]
\newtheorem{defn}{Definition}[theor]
\newtheorem{corollary}{Corollary}[section]
\newenvironment{proof}[1][Proof]{\begin{trivlist}
\item[\hskip \labelsep {\bfseries #1}]}{\end{trivlist}}

\centerline{\bf CONVERGENCE OF THE $J$-FLOW}
\centerline{\bf ON K\"AHLER SURFACES}
\bigskip
\centerline{Ben Weinkove}
\centerline{Department of Mathematics, Columbia University}
\centerline{New York, NY 10027}
\centerline{\it E-mail: weinkove@math.columbia.edu}
\bigskip

\setlength\arraycolsep{2pt}
\addtocounter{section}{1}

\bigskip
\noindent
{\bf 1. Introduction}
\bigskip

In \cite{Do}, Donaldson described how a number of geometric situations fit
into a general framework of diffeomorphism groups and moment maps.  In the
K\"ahler setting, he used this framework to define a natural parabolic
flow, as
follows. Suppose that
$(M, \omega)$ is a compact K\"ahler manifold of dimension
$n$ and let $\chi_0$ be another K\"ahler form on $M$, in a different
K\"ahler
class.  Consider the infinite-dimensional manifold $\mathcal{M}$ of
diffeomorphisms $f: M \rightarrow M$, homotopic to the identity.
$\mathcal{M}$
carries a natural symplectic form
$\Omega$ defined by
$$\Omega_f (v,w) = \int_M \omega(v,w) \frac{\chi_0^n}{n!},$$
for sections $v$, $w$ of $f^*(TM)$.
The group $\mathcal{G}$ of exact
$\chi_0$-symplectomorphisms of $M$ acts on $\mathcal{M}$ by
composition on the right, preserving $\Omega$.  We can identify the Lie
algebra
of $\mathcal{G}$ with the space of functions on $M$ of integral zero with
respect
to the volume form induced by $\chi_0$. A moment map $\mu: \mathcal{M}
\rightarrow \textrm{Lie}(\mathcal{G})^*$ for the group action is given by
$$\mu(f) = \frac{f^*(\omega) \wedge \chi_0^{n-1}}{\chi_0^n} - \frac{\int_M
\omega
\wedge \chi_0^{n-1}}{\int_M \chi_0^n},$$
where we are using the $L^2$ inner product to identify
$\textrm{Lie}(\mathcal{G})$ with its dual. It is natural to look for
solutions
of
\begin{equation} \label{eqnmoment}
\mu(f) = 0 \qquad (\textrm{mod} \ \mathcal{G}).
\end{equation}
These points form the symplectic quotient.  Under certain conditions, one
would
hope that the gradient flow $f_t$ of the function $\| \mu \|^2$ on
$\mathcal{M}$ would converge to give a solution of (\ref{eqnmoment}).
The
gradient flow can be rewritten as a flow of K\"ahler forms
 $(f_t^*)^{-1}(\chi_0)$ on $M$.  This defines a parabolic flow on
the space of K\"ahler potentials and is the object of study of this paper.

At around the same time, Chen \cite{C1} independently
discovered the same flow as the gradient flow of his
$J$-functional.  He later called it the $J$-flow \cite{C2}.
He showed in \cite{C1} 
that the
$J$-functional is related to the
Mabuchi K-energy
\cite{Ma}, which plays a key role in the study of K\"ahler geometry and
stability in the sense of geometric invariant theory (see \cite{Y2},
\cite{T2},
\cite{T3} and \cite{PS} for example).

Explicitly, the $J$-flow is defined as follows.
Let $c$ be the constant given by
$$c = \frac{ \int_M \omega \wedge \chi_0^{n-1}}{\int_M \chi_0^n},$$
and let $\mathcal{H}$ be the space of K\"ahler potentials
$$\mathcal{H} = \{ \phi \in C^{\infty}(M) \ | \ \chi_{\phi} = \chi_0 +
\frac{\sqrt{-1}}{2} \partial \dbar \phi >0 \}.$$
The $J$-flow is the flow on $\mathcal{H}$ given by
\begin{eqnarray}\nonumber
\ddt{\phi_t} & = & c - \frac{\omega \wedge
\chi_{\phi_t}^{n-1}}{\chi_{\phi_t}^n}.
\\
\label{eqnJflow}
\phi_0 & = & 0.
\end{eqnarray}
A critical point of the $J$-flow gives a K\"ahler metric $\chi$ satisfying
\begin{equation} \label{maineqn}
\omega \wedge \chi^{n-1} = c \chi^n.
\end{equation}
Donaldson \cite{Do} asked whether one can find a
solution to (\ref{maineqn}) in the class $[\chi_0]$ under certain
assumptions.  He
noted that a necessary condition is that $[nc\chi_0 - \omega]$ be a
K\"ahler
class, and conjectured that this condition be sufficient.  Chen \cite{C1}
confirmed this conjecture in the case $n=2$, without using the $J$-flow,
by
observing that (\ref{maineqn})  reduces to a Monge-Amp\`ere equation which
can be
solved by the well-known result of Yau \cite{Ya}.
The conjecture is still open for $n>2$.

Chen \cite{C1} shows that
Donaldson's conjecture would imply a result on the lower
bound of the Mabuchi K-energy for compact K\"ahler manifolds $M$ with
negative first Chern class.  Namely, if $-\omega \in c_1(M)$ with
$\omega>0$, then for K\"ahler classes $[\chi_0]$ satisfying
$$nc [\chi_0] - [\omega] >0,$$
the Mabuchi K-energy would have a lower bound in the class $[\chi_0]$.

Solutions of the $J$-flow exist for a short time by general
theory, since the flow is parabolic.  In \cite{C2}, Chen showed that the
flow
always exists for all time for any smooth initial data.  He also showed
that if
the bisectional curvature of
$\omega$ is non-negative then the $J$-flow converges to a critical metric.

In general, the behaviour of the flow is not known.
In this paper, we deal with the case $n=2$
with no curvature restrictions.  Our main result is as follows.

\bigskip
\noindent
{\bf Main Theorem} {\it Suppose that $(M, \omega)$ has dimension $n=2$ and
that
$$nc \chi_0 - \omega >0.$$
Then the $J$-flow (\ref{eqnJflow}) converges in $C^{\infty}$ to a smooth
critical
metric.}
\bigskip

The outline of the paper is as follows.  In section 2 we state some
preliminary
facts about the flow and introduce notation.  In section 3, the maximum
principle
is used to derive an estimate on the second derivatives of $\phi$ in
terms of $\phi$ itself.  In section 4, a $C^0$ estimate for $\phi$ is
given.
The argument uses the second order estimate, a Moser iteration argument
applied
to the exponential of $- \phi$ and the result of Tian \cite{Ti} (see also
\cite{TY}) on the existence of constants $\al>0$ and $C$ such that
$$\int_M e^{-\al \phi} \frac{\chi_0^n}{n!} \le C,$$
for all $\phi$ in $\mathcal{H}$ with $\sup_M \phi =0$.  In section 5,
the proof of the main theorem is completed.

\addtocounter{section}{1}
\setcounter{equation}{0}
\bigskip
\bigskip
\noindent
{\bf 2. Preliminaries and notation}
\bigskip

From now on, assume that $\omega$ has been scaled so that $c=1/n.$  We
will work
in local coordinates, and write
$$\omega = \frac{\sqrt{-1}}{2} g_{i \overline{j}} dz^i \wedge
dz^{\overline{j}},
\qquad
\chi_0 = \frac{\sqrt{-1}}{2} \chz{i}{j} dz^i \wedge dz^{\overline{j}},$$
and
$$\chi = \frac{\sqrt{-1}}{2} \ch{i}{j} dz^i \wedge dz^{\overline{j}}
=  \frac{\sqrt{-1}}{2} (\chz{i}{j} + \dd{i} \dbr{j} \phi) dz^i \wedge
dz^{\overline{j}},$$
where $\chi= \chi_{\phi}$ (suppressing the $t$-subscript.)  The
operators
$\Lambda_{\omega}$ and
$\Lambda_{\chi}$ act on $(1,1)$ forms $\al = \frac{\sqrt{-1}}{2}
\al_{i\overline{j}} dz^i
\wedge dz^{\overline{j}}$ by
$$ \Lambda_{\omega} \al = \gi{i}{j} \al_{i \overline{j}}, \qquad
\textrm{and}
\qquad \Lambda_{\chi} \al = \chii{i}{j} \al_{i \overline{j}}.$$

The $J$-flow (\ref{eqnJflow}) can be written
\begin{eqnarray} \nonumber
\ddt{\phi} & = & \frac{1}{n} (1- \Lambda_{\chi} \omega) \\
\label{eqnJflow2}
\phi|_{t=0} & = & 0.
\end{eqnarray}
Differentiating with respect to $t$ gives
\begin{eqnarray} \label{eqnddt}
\ddt{}{\left(\ddt{\phi}\right)} & = & \lt \left(\ddt{\phi}\right),
\end{eqnarray}
where the operator $\lt$ acts on functions $f$ by
$$\lt f = \frac{1}{n}\chii{k}{j}\chii{i}{l}\g{i}{j} \dd{k} \dbr{l}
f.$$
For convenience, write
$$\h{k}{l} = \chii{k}{j}\chii{i}{l}\g{i}{j}.$$
The tensor $\h{k}{l}$ is positive definite and its inverse
defines a Hermitian metric on $M$. The
operator
$\lt$ is, up to a constant factor, the Laplacian associated to
this Hermitian metric.

By the maximum principle for parabolic equations, (\ref{eqnddt}) implies
that
\begin{equation} \label{eqnupperlowerbound}
\inf_M (\Lambda_{\chi_0} \omega) \le \Lambda_{\chi}\omega \le \sup_M
(\Lambda_{\chi_0}
\omega),
\end{equation}
which gives a lower bound for $\chi$,
\begin{equation} \label{eqnlowerbound}
\chi \ge \frac{1}{\sup_M (\Lambda_{\chi_0} \omega)} \,
\omega.
\end{equation}

The $J$-functional \cite{C1} is defined by
$$J_{\omega, \chi_0} (\phi) = \int_0^1 \int_M \ddt{\phi_t} \,
\frac{\ \omega
\wedge \chi_{\phi_t}^{n-1}}{(n-1)!} \, dt,$$
where $ \{ \phi_t \}$ is a path in $\mathcal{H}$ between $0$ and $\phi$.
The functional is independent of the choice of path.
We will need the following formula for the functional in the case $n=2$.
Taking
the path $\phi_t = t \phi$, we see that
\begin{equation} \label{eqnJ}
J_{\omega, \chi_0}(\phi) = \frac{1}{2} \int_M \phi \, \omega \wedge
(\chi_0 +
\chi).
\end{equation}
Chen also makes use of the $I$-functional,
$$I_{\omega, \chi_0} (\phi) = \int_0^1 \int_M \ddt{\phi_t} \,
\frac{\chi_{\phi_t}^n}{n!} \, dt.$$
This is a well-known functional in K\"ahler geometry (see
\cite{Ma}). Notice that $I(\phi)=0$ along the flow.    For
$n=2$, this functional is given by
\begin{equation} \label{eqnI}
I_{\omega, \chi_0}(\phi) = \frac{1}{6} \int_M \phi \, ( \chi_0^2 + \chi
\wedge
\chi_0 + \chi^2).
\end{equation}

In the course of the paper, $C_0, C_1, \ldots$ will denote constants
depending
only on the initial data $\omega$ and $\chi_0$.  Curvature expressions
such as $\R{i}{j}{k}{l}$ will always refer to the metric $\g{i}{j}$.

\addtocounter{section}{1}
\setcounter{equation}{0}
\bigskip
\bigskip
\noindent
{\bf 3. Second order estimate}
\bigskip

We use the maximum principle to obtain an estimate on the second
derivative of
$\phi$ in terms of $\phi$.  We choose to calculate the evolution of $(\log
\Lambda_{\omega} \chi - A\phi)$ for some constant $A$ (compare to
\cite{Ya}, \cite{Au} or \cite{Si} for the analogous estimate for
the well-known Monge-Amp\`ere equation, and \cite{Ca} for the 
K\"ahler-Ricci flow.)

\begin{theorem} \label{theoremC2}
Suppose that $(M, \omega)$ has dimension $n=2$ and that
\begin{equation} \label{eqncondition}
\chi_0 - \omega >0.
\end{equation}
Let $\phi=\phi_t$ be a solution of the $J$-flow (\ref{eqnJflow2}) on $[0,
\infty)$.  Then there exist constants $A>0$ and $C>0$ depending only on
the
initial data such that for any time $t\ge0$, $\chi = \chi_{\phi_t}$
satisfies
\begin{equation} \label{eqnC2}
\Lambda_{\omega} \chi \le C e^{A (\phi - \inf_{M \times [0,t]}
\phi)}.
\end{equation}
\end{theorem}
\addtocounter{lemma}{1}

\begin{proof}
We will calculate
$$(\lt - \ddt{})(\log (\Lambda_{\omega} \chi) - A \phi).$$
Using normal coordinates for $\omega$, first calculate
\begin{eqnarray*}
\lt(\Lambda_{\omega} \chi) & = & \frac{1}{n} \h{k}{l} \dd{k} \dbr{l}
(\gi{i}{j} \ch{i}{j}) \\
& = & \frac{1}{n} \h{k}{l} R^{\ \ \, i \overline{j}}_{k \overline{l}}
\ch{i}{j} +
\frac{1}{n} \h{k}{l} \gi{i}{j}\dd{k}\dbr{l} \ch{i}{j}.
\end{eqnarray*}
And
\begin{eqnarray*}
\ddt{} (\Lambda_{\omega} \chi) & = & \ddt{} (\gi{i}{j} \dd{i}\dbr{j} \phi)
\\
& = & - \frac{1}{n} \gi{i}{j} \dd{i}\dbr{j} (\chii{k}{l} \g{k}{l}) \\
& = & \frac{1}{n} ( \gi{i}{j} \dd{i} (\chii{p}{l} \dbr{j}\ch{p}{q}
\chii{k}{q})\g{k}{l}
 +  \gi{i}{j} \chii{k}{l} \R{i}{j}{k}{l} )\\
& = & \frac{1}{n} ( \gi{i}{j} \h{p}{q} \dd{i}\dbr{j} \ch{p}{q} -
\gi{i}{j} \h{r}{q} \chii{p}{s} \dd{i} \ch{r}{s} \dbr{j}\ch{p}{q} \\
&& \mbox{} - \gi{i}{j} \h{p}{s} \chii{r}{q} \dd{i} \ch{r}{s}
\dbr{j}\ch{p}{q}
+ \chii{k}{l} \Ric{k}{l}).
\end{eqnarray*}
Now
$$\lt \log( \Lambda_{\omega} \chi) = \frac{\lt (\Lambda_{\omega}
\chi)}{\Lambda_{\omega} \chi} - \frac{ | \tilde{\nabla}
(\Lambda_{\omega} \chi)|^2}{(\Lambda_{\omega} \chi)^2},$$
where
$$| \tilde{\nabla} (\Lambda_{\omega} \chi)|^2 = \frac{1}{n} \h{k}{l}
\dd{k}
(\Lambda_{\omega} \chi) \dbr{l} (\Lambda_{\omega} \chi).$$
Note that by the K\"ahler property of $\chi$, we have
$$\dd{i}\dbr{j} \ch{k}{l} = \dd{k}\dbr{l} \ch{i}{j}.$$
Then
\begin{eqnarray*}
\lefteqn{(\lt - \ddt{}) \log(\Lambda_{\omega} \chi) } \\
& = & \frac{1}{n
\Lambda_{\omega}
\chi} ( \h{k}{l} R_{k \overline{l}}^{\ \ \, i \overline{j}} \ch{i}{j} - n
\frac{|
\tilde{\nabla} (\Lambda_{\omega} \chi)|^2}{\Lambda_{\omega} \chi}
+  \gi{i}{j} \h{r}{q} \chii{p}{s} \dd{i} \ch{r}{s} \dbr{j} \ch{p}{q} \\
&& \qquad \quad \mbox{} + \gi{i}{j} \h{p}{s} \chii{r}{q} \dd{i} \ch{r}{s}
\dbr{j}
\ch{p}{q} -  \chii{k}{l} \Ric{k}{l}).
\end{eqnarray*}

We need the following lemma to deal with the second term on the right hand
side.

\begin{lemma}
$$ n | \tilde{\nabla} (\Lambda_{\omega} \chi)|^2 \le
(\Lambda_{\omega} \chi) \gi{i}{j} \h{r}{q} \chii{p}{s} \dd{i} \ch{r}{s}
\dbr{j} \ch{p}{q}.$$
\end{lemma}
\begin{proof}
Using normal coordinates for $\omega$ in which $\chi$
is diagonal, and making use of the Cauchy-Schwartz inequality, we obtain
\begin{eqnarray*}
n | \tilde{\nabla} (\Lambda_{\omega} \chi)|^2 & = & \sum_{i,j,k}
\chii{k}{k}
\chii{k}{k}
\dd{k} \ch{i}{i} \dbr{k}\ch{j}{j} \\
& \le & \sum_{i,j} \left(\sum_{k} (\chii{k}{k})^2 |\dd{k}
\ch{i}{i}|^2\right)^{1/2} \left(\sum_{k} (\chii{k}{k})^2 |\dd{k}
\ch{j}{j}|^2\right)^{1/2}
\\ & = & \left(\sum_i \left(\sum_{k} (\chii{k}{k})^2 |\dd{k}
\ch{i}{i}|^2 \right)^{1/2}\right)^2
 \\ & = & \left(\sum_i \sqrt{\ch{i}{i}} \, \left(\sum_{k}
(\chii{k}{k})^2 \chii{i}{i} |\dd{k}
\ch{i}{i}|^2\right)^{1/2}\right)^2 \\
& \le & \sum_i \ch{i}{i} \sum_{i,k} (\chii{k}{k})^2 \chii{i}{i} |\dd{k}
\ch{i}{i}|^2 \\
& = & (\Lambda_{\omega} \chi) \sum_{i,k} (\chii{k}{k})^2 \chii{i}{i}
\dd{k}
\ch{i}{i} \dbr{k} \ch{i}{i} \\
& = & (\Lambda_{\omega} \chi) \sum_{i,k} (\chii{k}{k})^2 \chii{i}{i}
\dd{i}
\ch{k}{i} \dbr{i} \ch{i}{k} \\
& \le & (\Lambda_{\omega} \chi) \sum_{i,j,k} (\chii{k}{k})^2 \chii{i}{i}
\dd{j}
\ch{k}{i} \dbr{j} \ch{i}{k} \\
& = & (\Lambda_{\omega} \chi) \gi{i}{j} \h{r}{q} \chii{p}{s} \dd{i}
\ch{r}{s}
\dbr{j} \ch{p}{q}.
\end{eqnarray*}
\end{proof}
\bigskip

Let $C_0$ be a constant satisfying
$$R_{k \overline{l}}^{\ \ \, i \overline{j}} \ge - C_0 \g{k}{l}
\gi{i}{j}.$$
Then,
\begin{eqnarray*}
(\lt - \ddt{} ) \log (\Lambda_{\omega} \chi) & \ge & \frac{1}{n
\Lambda_{\omega} \chi} ( - C_0 \h{k}{l} \g{k}{l} \gi{i}{j} \ch{i}{j} -
\chii{k}{l} \Ric{k}{l})\\
& = & \frac{1}{n} ( -C_0 \h{k}{l} \g{k}{l} - \frac{1}{\Lambda_{\omega}
\chi}
\chii{k}{l} \Ric{k}{l}).
\end{eqnarray*}
Now calculate
\begin{eqnarray*}
(\lt - \ddt{} ) \phi & = & \frac{1}{n} ( \h{k}{l} \dd{k} \dbr{l} \phi +
\chii{i}{j} \g{i}{j} - 1) \\
& = & \frac{1}{n} (\chii{k}{j} \chii{i}{l} \g{i}{j} \ch{k}{l} - \h{k}{l}
\chz{k}{l} + \chii{i}{j} \g{i}{j} - 1) \\
& = & \frac{1}{n} (2 \chii{i}{j} \g{i}{j} - \h{k}{l} \chz{k}{l} - 1).
\end{eqnarray*}
At this point we must choose our value of $A$. From our
assumption (\ref{eqncondition}), we can choose
$0< \ep < 1/3$ to be sufficiently small so that
\begin{equation} \label{eqnchz}
\chi_0 \ge (1 + 3 \ep) \omega.
\end{equation}
Let $A$ be given by
$$A = \frac{C_0}{\ep}.$$

Fix a time $t>0$.  There is a point $(x_0, t_0)$ in $M \times [0,t]$ at
which
the maximum of $(\log (\Lambda_{\omega} \chi) - A\phi)$ is achieved.  We
may
assume that $t_0>0$.  At this point, we have
\begin{eqnarray*}
0 & \ge & (\lt - \ddt{} ) (\log(\Lambda_{\omega} \chi) - A \phi) \\
& \ge & \frac{1}{n} (- C_0 \h{k}{l} \g{k}{l} - \frac{1}{\Lambda_{\omega}
\chi} \chii{k}{l} \Ric{k}{l} - 2A \chii{i}{j} \g{i}{j} + A\h{k}{l}
\chz{k}{l} + A)
\\
& \ge & \frac{1}{n} (- C_0 \h{k}{l} \g{k}{l} - \frac{1}{\Lambda_{\omega}
\chi} \chii{k}{l} \Ric{k}{l} - 2A \chii{i}{j} \g{i}{j} + (1-\ep) A\h{k}{l}
\chz{k}{l} \\
&& \mbox{} + 
\ep A \h{k}{l} \g{k}{l} +  A)
\\
& = & \frac{1}{n} ( - \frac{1}{\Lambda_{\omega}
\chi} \chii{k}{l} \Ric{k}{l} - 2A \chii{i}{j} \g{i}{j} + (1-\ep) A\h{k}{l}
\chz{k}{l}  +  A).
\end{eqnarray*}
From the lower bound (\ref{eqnlowerbound}) on $\ch{k}{l}$, the term
$\chii{k}{l}
\Ric{k}{l}$ is bounded above and hence at $(x_0, t_0)$, we have
$$
1 + (1-\ep) \h{k}{l} \chz{k}{l} - 2\chii{i}{j} \g{i}{j} \le
\frac{C_1}{ (\Lambda_{\omega} \chi)}.
$$
From (\ref{eqnchz}), we get
\begin{equation} \label{eqnatmax}
1+ (1+\ep) \h{k}{l} \g{k}{l} - 2\chii{i}{j} \g{i}{j} \le \frac{C_1}{
(\Lambda_{\omega} \chi)}.
\end{equation}
 We will
compute in normal coordinates at
$x_0$ for
$\omega$ in which $\chi$ is diagonal and has eigenvalues $\lambda_1,
\lambda_2$.  From (\ref{eqnlowerbound}), $\lambda_1$ and $\lambda_2$ are
bounded
below by a positive constant.  We want to show that they are also bounded
above.
First, observe that for $n=2$,
$$ \frac{1}{\Lambda_{\chi} \omega} = \frac{\det \chi}{ (\det{\omega})
(\Lambda_{\omega}
\chi)},$$
and by (\ref{eqnupperlowerbound}), this is bounded along the flow.

Multiplying (\ref{eqnatmax}) by $(\det \chi / \det \omega)$
gives,
$$ \lambda_1 \lambda_2 + (1+ \ep)( \frac{\lambda_2}{\lambda_1} +
\frac{\lambda_1}{\lambda_2})
 - 2(\lambda_1 + \lambda_2) \le C_2.$$
From (\ref{eqnupperlowerbound}), we may suppose that one of the
eigenvalues, say
$\lambda_2$,  is bounded from above.  Rewrite the
inequality as
$$\lambda_1 (\lambda_2 + (1+ \ep)\frac{1}{\lambda_2} -2) + (1+ \ep)
\frac{\lambda_2}{\lambda_1} - 2\lambda_2 \le C_2.$$
Then, since the function $f : \, (0,\infty) \rightarrow \mathbf{R}$
defined by
$$f(x) = x+ (1 + \ep)\frac{1}{x} - 2,$$
is bounded below by a small positive constant depending on $\ep$, we see
that
$\lambda_1$ must also be bounded above.  Hence at the point $(x_0, t_0)$,
there
exists $C$ depending only on the initial data such that
$$\Lambda_{\omega} \chi \le C.$$
Then, on $M \times [0,t]$,
$$\log (\Lambda_{\omega} \chi) - A\phi \le \log C - A\inf_{M \times [0,t]}
\phi.$$
Exponentiating gives
$$\Lambda_{\omega} \chi \le C e^{A(\phi - \inf_{M \times [0,t]} \phi)},$$
completing the proof of the theorem.
\end{proof}

\addtocounter{section}{1}
\setcounter{theorem}{0}
\setcounter{equation}{0}
\setcounter{lemma}{0}
\bigskip
\pagebreak[3]
\noindent
{\bf 4. Zero order estimate}
\bigskip

We prove an estimate on the $C^0$ norm of $\phi$ using a Moser iteration
method
applied to the exponential of the solution rather than a power of the
solution
(compare to \cite{Ya}) and the estimate of Theorem
\ref{theoremC2}.

\addtocounter{lemma}{1}
\addtocounter{proposition}{1}
\begin{theorem} \label{theoremC0}
Suppose that $(M, \omega)$ has dimension $n=2$ and that
$$\chi_0 - \omega >0.$$
Let $\phi_t$ be a solution of the $J$-flow (\ref{eqnJflow2}) on $[0,
\infty)$.
Then there exists
a constant $\tilde{C}$ depending only on the initial data such that
$$\| \phi_{t} \|_{C^0(M)} \le \tilde{C}.$$
\end{theorem}
\begin{proof}
Suppose first that $\inf_M \phi_t$ is bounded from below
uniformly in time.  We will show that this implies the above estimate.
Since
the functional $J_{\omega, \chi_0}$ decreases along the flow, there exists
a constant $C_0$ such that
$$ \int_M \phi_t \, \omega \wedge (\chi_0 + \chi_{\phi_t}) \le C_0,$$
using (\ref{eqnJ}).  Let $C_1$ be a positive constant satisfying
$$\omega^2 \le C_1 \omega \wedge \chi_0.$$
Then 
\begin{eqnarray*}
\int_M \phi_t \, \omega^2 & = & \int_M (\phi_t -\inf_M \phi_t) \omega^2 + \int_M
\inf_M \phi_t \, \omega^2 \\
& \le & C_1 \int_M (\phi_t - \inf_M \phi_t ) \omega \wedge \chi_0 + \inf_M
\phi_t  \int_M
\omega^2 \\
& \le & C_1 C_0 - C_1 \int_M \phi_t \, \omega \wedge \chi_{\phi_t} + \inf_M \phi_t
\left( \int_M \omega^2 - C_1 \int_M \omega \wedge \chi_0 \right) \\
& = & C_1 C_0 - C_1\int_M (\phi_t -\inf_M \phi_t) \omega \wedge \chi_{\phi_t} \\
&& \mbox{} 
+
\inf_M \phi_t \left(\int_M \omega^2 - 2C_1 \int_M \omega \wedge \chi_0 \right) \\
& \le & C_1 C_0 +
\inf_M \phi_t \left(\int_M \omega^2 - 2C_1 \int_M \omega \wedge \chi_0 \right).
\end{eqnarray*}
This gives an upper bound for $\int_M
\phi_t \, \omega^2$ depending on the lower bound for $\inf_M \phi_t$. 
Since $\triangle_{\omega} \phi_t > - \Lambda_{\omega} \chi_0$ along the
flow, it follows from the existence of a lower bound on the Green's
function of
$\omega$ that $\sup_M \phi_t$ is bounded from above, giving us the
required
estimate.

Now suppose that no such lower bound for $\inf_M \phi_t$ exists.  Then
 we can assume that there is a sequence of times $t_i
\rightarrow
\infty$ such that
\begin{enumerate}
\item[(i)] $\inf_M \phi_{t_i} = \inf_{t \in [0,t_i]} \inf_M \phi_t$
\item[(ii)] $\inf_M \phi_{t_i} \rightarrow - \infty$.
\end{enumerate}
We will seek a contradiction.
For a
fixed $i$, write
$$\psi_{t_i} = \phi_{t_i} - \sup_M \phi_{t_i}.$$
Notice that  $\sup_M \phi_{t_i}$ is bounded from below by zero from
$(\ref{eqnI})$
and the fact that $I(\phi_t)=0$.  Hence
$$\| \psi_{t_i} \|_{C^0} \rightarrow  \infty.$$
The following proposition is the key result of this section.  

\begin{proposition} \label{keyprop} Let $M$ be a compact complex surface with two
K\"ahler metrics $\chi_0$ and $\omega$.  Suppose that $\psi \in C^{\infty}(M)$
satisfies the conditions
$$\chi_{\psi} = \chi_0 + \frac{\sqrt{-1}}{2} \partial \dbar \psi >0, \qquad \sup_M
\psi =0,$$
and
$$\Lambda_{\omega} \chi_{\psi} \le C e^{A(\psi - \inf_M \psi)}.$$
Then there exists a constant $C'$ depending
only on $M$, $\omega$, $\chi_0$ and the constants $A$ and $C$ such that
$$\| \psi \|_{C^0} \le C'.$$
\end{proposition}

We apply this proposition to $\psi = \psi_{t_i}$ and obtain a
contradiction since 
\begin{eqnarray*}
\Lambda_{\omega} \chi_{\psi_{t_i}} & = & \Lambda_{\omega} \chi_{\phi_{t_i}} \\
& \le & C e^{A(\phi_{t_i} - \inf_{t \in [0,t_i]} \inf_M \phi_t)} \\
& = & Ce^{A(\psi_{t_i} - \inf_M \psi_{t_i})},
\end{eqnarray*}
where we have used 
Theorem \ref{theoremC2} and condition (i) above.  It remains to prove the
proposition.
\end{proof}

\bigskip
\noindent
{\bf Proof of Proposition \ref{keyprop}} \ Let $\delta$ be a small positive 
constant, to
be determined later. Set
$B = A/(1-\delta)$ and 
let $u=e^{-B\psi}$.  

Now, for $\be = n/(n-1) =2$, the
Sobolev inequality for functions $f$ on $(M, \omega)$ is
$$ \| f \|^2_{2\be} \le C_2 (\| \nabla f \|^2_{2} + \| f \|^2_{2}),$$
for $C_2$ depending on $\omega$.  We will apply this to $u^{p/2}$ for
$p \ge 1$.  This gives
\begin{equation} \label{eqnsobolev}
\left( \int_M e^{-Bp\be \psi} \frac{\omega^2}{2} \right)^{1/\be} \le C_2
\left(\int_M |
\nabla e^{-Bp\psi/2} |^2 \frac{\omega^2}{2} + \int_M e^{-Bp\psi}
\frac{\omega^2}{2} \right).
\end{equation}
Now calculate
\begin{eqnarray*}
\int_M | \nabla e^{-Bp\psi/2}|^2 \frac{\omega^2}{2} & = & \sqrt{-1} \int_M
\partial e^{-Bp\psi/2}
\wedge \dbar e^{-Bp\psi/2} \wedge \omega \\
& = & \frac{B^2 p^2}{4} \sqrt{-1} \int_M e^{-Bp\psi} \partial \psi \wedge
\dbar
\psi \wedge \omega \\
& = & - \frac{Bp}{4} \sqrt{-1} \int_M \partial (e^{-Bp\psi}) \wedge \dbar
\psi
\wedge \omega \\
& = & \frac{Bp}{2} \int_M e^{-Bp\psi} \frac{\sqrt{-1}}{2} \partial \dbar
\psi \wedge \omega \\
& = & \frac{Bp}{2} \int_M e^{-Bp\psi} (\chi_{\psi} - \chi_0) \wedge \omega \\
& = & \frac{Bp}{2} \int_M e^{-Bp \psi} (\Lambda_{\omega} \chi_{\psi} -
\Lambda_{\omega} \chi_0) \frac{\omega^2}{2} \\
& \le & \frac{CBp}{2} \int_M e^{-Bp\psi} e^{A(\psi - \inf_M \psi)}
\frac{\omega^2}{2}
\\ & = & \frac{CBp}{2} e^{-A\inf_M{\psi}} \int_M e^{-(p-(1-\delta))B\psi}
\frac{\omega^2}{2},
\end{eqnarray*}
where we have used the estimate
$$\Lambda_{\omega} \chi_{\psi} \le 
C
e^{A(\psi - \inf_M \psi)}.$$

Then in (\ref{eqnsobolev}),
$$ \left(\int_M u^{p\be} \frac{\omega^2}{2} \right)^{1/\be} \le C_3 p
e^{-A\inf_M
\psi}
\int_M u^{p - (1-\delta)}
\frac{\omega^2}{2}.$$
Raising to the power $1/p$ and writing $\ga = 1-\de$ gives
$$ \| u \|_{p\be} \le C_3^{1/p} p^{1/p} e^{-(A/p)\inf_M \psi} \| u \|_{p -
\ga}^{(p - \ga)/p}.$$
Take the logarithm of both sides to get
$$\log \|u \|_{p\be} \le \frac{1}{p} \log C_3 + \frac{1}{p} \log p +
\frac{1}{p}
\sup_M(-A\psi) + \frac{(p- \ga)}{p} \log\|u \|_{p-\ga}.$$
We now apply the iteration.  First, replace $p$ with $p \be + \ga$ to get
\begin{eqnarray*}
\log \| u \|_{p\be^2 + \ga \be} & \le & \frac{1 + \be}{p\be + \ga} \log
C_3 +
\frac{1}{p\be + \ga}(\be \log p + \log(p\be + \ga)) \\
&& \mbox{} + \frac{1 + \be}{p\be + \ga} \sup_M(-A \psi) + \frac{\be
(p-\ga)}{p\be +
\ga} \log \|u\|_{p-\ga}.
\end{eqnarray*}
Repeat this procedure, replacing $p$ with $p\be + \ga$ to obtain for any
positive
integer $k$,
\begin{eqnarray}
\nonumber
\lefteqn{ \log \|u \|_{ p \be^{k+1} + \ga(\be + \be^2 + \ldots + \be^k)}}
\\
\nonumber
& \le &  \frac{1+
\be + \be^2 + \ldots + \be^k}{p\be^k + \ga(1 + \be +\be^2 + \ldots +
\be^{k-1})}
\log C_3 \\
\nonumber
&& \mbox{} + \frac{1}{p \be^k + \ga(1 + \be + \ldots + \be^{k-1})}
\left ( \right. \be^k
\log p + \be^{k-1} \log (p \be + \ga) + \ldots \\
\nonumber
&& \mbox{} \qquad \ldots + \log (p\be^k + \ga(1+ \be + \ldots
+ \be^{k-1})  \left. \right ) \\
\nonumber
&& \mbox{} + \frac{1+
\be + \be^2 + \ldots + \be^k}{p\be^k + \ga(1 + \be +\be^2 + \ldots +
\be^{k-1})}
\sup_M (-A\psi) \\
&& \mbox{} + \frac{\be^k (p-\ga)}{p\be^k + \ga(1 + \be +\be^2 + \ldots +
\be^{k-1})} \log\| u \|_{p-\ga}. \label{eqnmoser}
\end{eqnarray}
Now set $p = 1+\de$.  Then, since $\be =2$ we have
$$ p\be^k + \ga(1 + \be + \be^2 + \ldots + \be^{k-1}) = 1 + \be + \be^2 +
\ldots
+ \be^k + \de.$$
Notice that the second term on the right hand side of (\ref{eqnmoser}) is
bounded
by
\begin{eqnarray*}
\log p + \frac{1}{\be} \log{\be^2}  + \ldots 
+
\frac{1}{\be^{k}} \log(\be^{k+1}) & \le & \log p + \log \be ( \sum_{i=1}^k
\frac{i+1}{\be^i})\\
& \le & C_4.
\end{eqnarray*}
Then
\begin{eqnarray*}
\lefteqn{\log \|u \|_{ p \be^{k+1} + \ga(\be + \be^2 + \ldots + \be^k)} }
\\ && \qquad \qquad \le \log C_3 + C_4 + \sup_M(-A\psi) + 2\de \max( \log
\|u\|_{2\de},0).
\end{eqnarray*}
Using the fact that $A= (1-\de)B$ and $-B\psi = \log u$, and letting
$k$ tend to infinity,
$$\log \|u \|_{C_0} \le C_5 + 2 \max( \log \|u\|_{2\de},0).$$
Hence we get the following inequality for $\psi$,
\begin{equation} \label{eqnC0}
\| \psi \|_{C^0} \le C_6 + C_7 \max\left( \log\left(\int_M e^{-2\de
B
\psi}
\frac{\omega^2}{2} \right)^{1/2\de},0 \right).
\end{equation}

We can now finish the estimate.  First, define
$$P(M, \chi_0) = \{ \Phi \in C^2(M) \ | \ \chi_0 +
\frac{\sqrt{-1}}{2} \partial \dbar
\Phi \ge 0, \ \sup_M{\Phi}=0 \}.$$
Then Proposition 2.1 of \cite{Ti} (see section 4.4, \cite{Ho}) states that
there exist constants
$\al>0$ and
$C_8$ depending only on $(M,\chi_0)$ such that
$$ \int_M e^{-\al \Phi} \frac{\chi_0^n}{n!} \le C_8 \qquad \textrm{for all
} \Phi
\in P(M,
\chi_0).$$
Define $\de$ to be
$$\de = \min \{ \frac{\al}{4A}, \frac{1}{2} \} >0.$$
Then the required estimate follows from (\ref{eqnC0}), since $\psi$ belongs to
$P(M, \chi_0)$.

\addtocounter{section}{1}
\setcounter{theorem}{0}
\setcounter{lemma}{0}
\setcounter{equation}{0}
\bigskip
\bigskip
\pagebreak[3]
\noindent
{\bf 5. Convergence of the flow}
\bigskip

In this section we complete the proof of the main theorem.  We assume,
using
the result of \cite{C2}, that a solution $\phi = \phi_t$ for the $J$-flow
exists for all time.  From Theorem
\ref{theoremC2} and Theorem \ref{theoremC0} we have uniform estimates on
$\phi$
and the derivatives $\dd{i} \dbr{j} \phi$, using the fact that
$$\ch{i}{j} = \chz{i}{j} + \dd{i}\dbr{j} \phi >0.$$
Since the operator
$$\frac{1}{n} (1 - \Lambda_{\chi} \omega),$$
is concave in the $\ch{i}{j}$,
it is well known that, by the work of Evans
\cite{E1, E2} and Krylov \cite{Kr} (see also \cite{Tr}), one can deduce a
uniform
H\"older estimate on the second derivatives
$\dd{i} \dbr{j} \phi.$
By differentiating the equation (\ref{eqnJflow2}) and applying standard
Schauder estimates for parabolic equations (see \cite{LSU} for example),
one can
obtain uniform estimates on all of the derivatives of $\phi$.  It then
follows
that there is a sequence of times $t_j \rightarrow \infty$ such that
$\phi_{t_j}$
converges in $C^{\infty}$ to some smooth function $\phi_{\infty}$.  In
order to
show that we have convergence without having to pass to a subsequence, we
will
use a modification of the argument in \cite{Ca}.

Notice that $\partial \phi/ \partial t$ satisfies the heat equation
$$\ddt{}{\left(\ddt{\phi}\right)} = \lt \left(\ddt{\phi}\right).$$
Since we have uniform bounds for $\ch{i}{j}$ from above and away from
zero, and
bounds on $\ddt{} \ch{i}{j}$ and all the covariant derivatives
of $\ch{i}{j}$ and $\ddt{} \ch{i}{j}$, it follows from  the Harnack
inequality of
Li and Yau \cite{LY} and the
argument in \cite{Ca} that
there exist positive constants $C_0$ and $\eta$, which are independent of
$t$,
such that
$$\sup_{M} \left( \ddt{\phi} \right) - \inf_{M} \left( \ddt{\phi}
\right) \le C_0 e^{-\eta t}.$$ Since
$$\int_M \ddt{\phi} \chi^2 =0,$$
$\partial \phi / \partial t$ must take on the value zero somewhere on $M$
for each
$t$, and so
$$\left| \ddt{\phi} \right| \le C_0 e^{-\eta t}.$$
Hence for any $0 < s < s'$, and any $x \in M$,
\begin{eqnarray*}
| \phi (x, s') - \phi (x,s)| & = & | \int_s^{s'} \ddt{\phi} (x,t) dt| \\
& \le & \int_s^{s'} | \ddt{\phi}(x,t) | dt \\
& \le & C_0 \int_s^{s'} e^{-\eta t} dt \\
&  =& C_0 \frac{1}{\eta} (e^{-\eta s} - e^{- \eta s'}),
\end{eqnarray*}
which tends to zero as $s$ and $s'$ tend to infinity.  Hence $\phi_t$
converges
in the $C_0$ norm to $\phi_{\infty}$.  It must converge also in the
$C^{\infty}$
topology, since otherwise there would exist an integer $N$, an
$\epsilon>0$ and a sequence $t_j \rightarrow \infty$ with
$$ \| \phi_{t_j} - \phi_{\infty} \|_{C^N} \ge \epsilon.$$
Since $\phi$ is bounded in all the $C^k$ norms, one could pass to a
subsequence of the $\phi_{t_j}$ which would converge to some
$\phi_{\infty}' \neq
\phi_{\infty}$, giving the contradiction.  This completes the proof.

\bigskip
\noindent
{\bf Acknowledgements.}  This work was completed while the author was a graduate
student at Columbia University, and these results form part of his PhD thesis
\cite{We}.  The author is very grateful to his advisor D.H.
Phong for his constant support and advice.  He also thanks  
 Jacob
Sturm and Jian Song for some helpful conversations, and the referee for some
constructive comments.


\small
\begin{thebibliography}{99}
\bibitem[Au]{Au} Aubin, T. {\em Equations du type Monge-Amp\`ere sur les
vari\'et\'es K\"ahleriennes compacts}, Bull. Sc. Math. 102 (1978), 119-121
\bibitem[Ca]{Ca} Cao, H-D. {\em Deformation of K\"ahler metrics to
K\"ahler-Einstein metrics on compact K\"ahler manifolds}, Invent. Math. 81
(1985), 359-372
\bibitem[C1]{C1} Chen, X. X. {\em On the lower bound of the Mabuchi energy
and
its application}, Int. Math. Res. Notices 12 (2000), 607-623
\bibitem[C2]{C2} Chen, X. X. {\em A new parabolic flow in K\"ahler
manifolds},
preprint, arXiv: math.DG/0009247
\bibitem[Do]{Do} Donaldson, S. K. {\em Moment maps and diffeomorphisms},
Asian J. Math. 3, No. 1 (1999), 1-16
\bibitem[E1]{E1} Evans, L. C. {\em Classical solutions of fully nonlinear,
convex,
second order elliptic equations}, Comm. Pure Appl. Math. 25 (1982),
333-363
\bibitem[E2]{E2} Evans, L. C. {\em Classical solutions of the
Hamilton-Jacobi
Bellman equation for uniformly elliptic operators}, Trans. Amer. Math.
Soc. 275
(1983), 245-255
\bibitem[Ho]{Ho} H\"ormander, L. {\em An introduction to complex analysis
in several variables}, Van Nostrand, Princeton, NJ 1973
\bibitem[Kr]{Kr} Krylov, N. V. {\em Boundedly nonhomogeneous elliptic and
parabolic equations}, Izvestia Akad. Nauk. SSSR 46 (1982), 487-523.
English
translation in Math. USSR Izv. 20 (1983), No. 3, 459-492
\bibitem[LSU]{LSU} Ladyzenskaja, O. A., Solonnikov, V. A. and Ural'Ceva,
N. N.
{\em Linear and quasilinear equations of parabolic type}, Providence,
Amer. Math.
Soc. 1968
\bibitem[LY]{LY} Li, P. and Yau, S.-T. {\em On the parabolic kernel of the
Schr\"odinger operator}, Acta Math. 156 (1986), No. 3-4, 153-201
\bibitem[Ma]{Ma} Mabuchi, T. {\em K-energy maps integrating Futaki
invariants},
T\^{o}hoku Math. Journ., 38 (1986), 575-593
\bibitem[PS]{PS} Phong, D. H. and Sturm, J. {\em Stability, energy
functionals,
and K\"ahler-Einstein metrics}, Comm. Anal. Geom. 11 (2003), No. 3, 565-597
\bibitem[Si]{Si} Siu, Y.-T. {\em Lectures on Hermitian-Einstein metrics
for
stable  bundles and K\"{a}hler-Einstein metrics}, Birkh\"{a}user Verlag,
Basel 1987
\bibitem[T1]{Ti} Tian, G. {\em On K\"ahler-Einstein metrics on certain
K\"ahler
manifolds with $c_1(M)>0$}, Invent. Math. 89 (1987), 225-246
\bibitem[T2]{T2} Tian, G. {\em The K-energy on hypersurfaces and
stability},
Comm. Anal. Geom. 2 (1994), No. 2, 239-265
\bibitem[T3]{T3} Tian, G. {\em K\"ahler-Einstein metrics with positive
scalar
curvature}, Invent. math. 137 (1997), 1-37
\bibitem[TY]{TY} Tian, G. and Yau, S.-T. {\em K\"ahler-Einstein metrics on
complex surfaces with $c_1(M)$ positive}, Comm. Math. Phys. 112 (1987),
\bibitem[Tr]{Tr} Trudinger, N. S. {\em Fully nonlinear, uniformly elliptic
equations under natural structure conditions}, Trans. Amer. Math. Soc. 278
(1983), 751-769
\bibitem[We]{We} Weinkove, B. {\em The J-flow, the Mabuchi energy, the Yang-Mills
flow and multiplier ideal sheaves}, PhD thesis, Columbia University 2004
\bibitem[Y1]{Ya} Yau, S.-T. {\em On the Ricci curvature of a compact
K\"ahler
manifold and the complex Monge-Amp\`ere equation, I}, Comm. Pure Appl.
Math. 31
(1978), 339-411
\bibitem[Y2]{Y2} Yau, S.-T. {\em Open problems in geometry}, Proc.
Symposia Pure
Math. 54 (1993), 1-28 (problem 65)
\end{thebibliography}
\end{document}